\documentclass{ps}
\usepackage{amsmath,amsfonts}
\usepackage{color}
\usepackage[applemac]{inputenc}
\usepackage{mathrsfs}

\def \O{\mathcal{O}}

\def\scD{{\mathscr D}}

\def\leftB{[\![}
\def\rightB{]\!]}

\newcommand \A[1]{{\bf (#1)}}

\def\O{{\cal{O}}}

\def\F{{\cal F}}

\def\P{{\mathbb{P}}  }

\def\bint#1^#2{\displaystyle{\int_{#1}^{#2}}}
\def\bsum#1^#2{\displaystyle{\sum_{#1}^{#2}}}

\def\xdt_#1{X_#1(\Delta t)}

\newtheorem{THM}{Theorem}[section]

\newtheorem{PROP}{Proposition}[section]

\newtheorem{REM}{Remark}[section]

\newcommand{\mysection}{\setcounter{equation}{0} \section}

\newcommand{\blackboard}[1]{\mathbf{#1}} 
\newcommand{\R}{\blackboard{R}}%		reels
\newcommand{\N}{\blackboard{N}}%		entiers
\renewcommand{\P}{\blackboard{P}}%	 	probabilite
\newcommand{\E}{\blackboard{E}}% 		esperance
\newcommand{\delimleft}[2]{\ifcase #1\or%	delimiteur gauche
    \bigl#2\or %
    \Bigl#2\or %
    \biggl#2\or %
    \Biggl#2\or %
    \left#2\fi}
\newcommand{\delimright}[2]{\ifcase #1\or%	delimiteur droite
    \bigr#2\or %
    \Bigr#2\or %
    \biggr#2\or %
    \Biggr#2\or %
    \right#2\fi}

\begin{document}
%%-----------------------------
%%      the top matter
%%-----------------------------
\title{Concentration Bounds for Stochastic Approximations}% At most 5 thanks
\author{N. Frikha}\address{LPMA, Universit\'e Denis Diderot, 175 Rue du Chevaleret 75013 Paris, frikha@math.univ-paris-diderot.fr}
\author{S. Menozzi}\address{Laboratoire d'Analyse et Probabilit\'es, Universit\'e d'Evry Val d'Essonne, 23 Boulevard de France, 91037 Evry Cedex,\\ menozzi@cmap.polytechnique.fr}

\date{\today}
\begin{abstract} We obtain non asymptotic concentration bounds for two kinds of stochastic approximations. 
We first consider the deviations 
between the expectation of a given function of an Euler like discretization scheme of some diffusion process at a fixed deterministic time and its empirical mean obtained by the Monte-Carlo procedure. We then give some estimates concerning the deviation between the value at a given time-step of a stochastic approximation algorithm and its target. Under suitable assumptions both concentration bounds turn out to be Gaussian.
 The key tool consists in exploiting accurately the concentration properties of the increments of the schemes. Also, no specific non-degeneracy conditions are assumed.
\end{abstract}
\subjclass{60H35,65C30,65C05}
\keywords{Non asymptotic bounds, Euler scheme, Stochastic approximation algorithms, Gaussian concentration}
\maketitle

\section{Statement of the Problem}
Let us consider a $d$-dimensional stochastic evolution scheme of the form\begin{equation}
\label{SCHEME_GEN}
\xi_{n+1}=\xi_n+\gamma_{n+1} F(n,\xi_n,Y_{n+1}), \ n\ge 0, \xi_0=x\in \R^d,
\end{equation}

\noindent where  $(\gamma_n)_{n\ge 1} $ is a deterministic positive sequence of time steps, the function $F:\N\times\R^d\times \R^q\rightarrow \R^d $ is a measurable function satisfying some assumptions that will be specified later on, and the $(Y_i)_{i\in \N^*}$ are i.i.d. $\R^q $-valued random variables defined on some probability space $(\Omega,\F,\P) $ whose law satisfies a \textit{Gaussian concentration property}. That is, there exists $\alpha>0$ s.t. for every real-valued 1-Lipschitz function $f$ defined on $\R^q$ and for all $\lambda\ge 0$:
\begin{equation}
\label{CONC_GAUSS}
 \E[\exp(\lambda f(Y_1))]\le \exp(\lambda \E[f(Y_1)]+\frac{\alpha\lambda^2}{4})\tag{$GC(\alpha)$}.
\end{equation}
From the Markov exponential inequality and \eqref{CONC_GAUSS}, one derives ${\cal D}(f,r):=\P[f(Y_1)-\E[f(Y_1)]\ge r]\le \exp(-\lambda r+\frac{\alpha\lambda^2}{4}), \forall \lambda ,r\ge 0 $. An optimization over $\lambda $ gives that ${\cal D}(f,r) $ %the deviation above $r$ of $f$ from its mean under the measure $\mu$ 
 has sub-Gaussian tails bounded by $\exp(-\frac{r^2}{\alpha}) $.

A practical criterion for \eqref{CONC_GAUSS} to hold is given by Bolley and Villani \cite{boll:vill:05}. If there exists $\varepsilon>0$ s.t. 
$\E[\exp(\varepsilon |Y_1|^2)] <+\infty$, then the law of $Y_1$ satisfies \eqref{CONC_GAUSS}
with $\alpha:=\alpha(\varepsilon)$. The two claims are actually equivalent.
In the following \eqref{CONC_GAUSS} is the only crucial property we require on the innovations $(Y_i)_{i\in \N^*} $. In particular we do not assume any absolute continuity of the law of $Y_1$ w.r.t. the Lebesgue measure.

We are interested in giving non asymptotic concentration bounds for two specific problems related to evolutions of type \eqref{SCHEME_GEN}. We first want to control the deviations of the empirical mean associated to a function of an Euler like discretization scheme of a diffusion process at a fixed deterministic time from the real mean. Secondly, we want to derive deviation estimates between the value of a Robbins-Monro type stochastic algorithm taken at fixed time-step and its target. Under some mild assumptions, we show that the Gaussian concentration property %\eqref{CONC_GAUSS}
 of the innovations transfers to the scheme. %The concentration of the averaged algorithm of Ruppert and Polyak is also investigated.
Concerning stochastic algorithms, our deviation results are to our best knowledge the first of this nature.
\subsection{Euler like Scheme of a Diffusion Process}
\label{DESCR_EUL}
Let $(\Omega,\F,(\F_t)_{t\ge 0},\P) $ be a filtered probability space satisfying the usual conditions and $(W_t)_{t\ge 0} $ be a $q$-dimensional $(\F_t)_{t\ge 0} $ Brownian motion. Let us  consider a $d$-dimensional diffusion process $(X_t)_{t\ge 0} $ with dynamics:
\begin{equation}
\label{EqDfSt}
X_t=x+\int_0^t b(s,X_s)ds+\int_0^t\sigma(s,X_s)dW_s,
\end{equation}
where the coefficients $b,\sigma $ are assumed to be uniformly Lipschitz continuous in space and measurable in time.

%``smooth" and to satisfy a non-degeneracy condition. Precisely, we suppose the following conditions are in force
%\begin{trivlist}
%\item[\A{S}] The coefficients $b:\R^+\times \R^d\rightarrow \R^d$ and $\sigma:\R^+\times \R^d\rightarrow \R^d\otimes \R^d$ are bounded and Lipschitz continuous in space. %twice continuously differentiable in space, measurable in time. 
%\item[\A{UE}] The diffusion coefficient is uniformly elliptic, i.e. there exists $\Lambda\in (0,1] $ s.t. forall $(t,x,\xi)\in \R^+\times (\R^d)^2$, $ \langle a(t,x)\xi,\xi\rangle \ge \Lambda |\xi|^2$.
%\end{trivlist}
%From now on we say that assumption \A{A} is satisfied if \A{S} and \A{UE} hold. 

For a given Lipschitz continuous function $f$ and a fixed deterministic time horizon $T$, quantities like $\E_x[f(X_T)] $ appear in many applications. In mathematical finance, it represents the price of a European option with maturity $T$ when the dynamics of the underlying asset is given by \eqref{EqDfSt}. Under suitable assumptions on the function $f$ and the coefficients $b,\sigma $, namely smoothness or non degeneracy, it can also be related to the Feynman-Kac representation of the heat equation associated to the generator of $X$. Two steps are needed to approximate $\E_x[f(X_T)] $:
\begin{trivlist}
\item[-] The first step consists in approximating the dynamics by a discretization scheme that can be simulated. 
For a given time step $\Delta=T/N,\ N\in \N^* $, setting for all $i\in\N,\ t_i:=i\Delta$,
we consider an Euler like scheme of the form:
\begin{eqnarray}
\label{EULER}
X_0^\Delta=x,\
\forall i\in \leftB 0, N-1\rightB, X_{t_{i+1}}^\Delta=X_{t_i}^\Delta+ b(t_i,X_{t_i}^\Delta)\Delta +\sigma(t_i,X_{t_i}^\Delta)\sqrt \Delta Y_{i+1},
\end{eqnarray}
\noindent where the $(Y_i)_{i\in \N^*} $ are $\R^q $-valued i.i.d. random variables whose law satisfies \eqref{CONC_GAUSS} for some $\alpha>0 $.  We also assume $\E[Y_1]=0_q,\ \E[Y_{1}Y_1^*]=I_{q} $, where $Y_1^* $ stands for the transpose of the column vector $Y_1$ and  $0_q,I_q $ respectively stand for the zero vector of $\R^q$ and the identity matrix of $\R^q\otimes \R^q$. 
%%%% A garder pour la suite
The previous assumptions include the case of the \textit{standard} Euler scheme, corresponding to $Y_1\overset{{\rm law}}{=}{\cal N}(0,I_q) $, which yields \eqref{CONC_GAUSS} with $\alpha=2 $ and the Bernoulli law $Y_1\overset{{\rm law}}{=} (B_1,\cdots, B_q),\ (B_k)_{k\in \leftB 1,q\rightB}$ i.i.d with law $\mu=\frac12 (\delta_{-1}+\delta_1)$. This latter choice can turn out to be useful, in terms of computational effort, to approximate \eqref{EqDfSt} when the dimension is large. 
\item[-] The second step consists in approximating the expectation $\E_x[f(X_T^\Delta)]$ involving the scheme \eqref{EULER} by a Monte-Carlo estimator:
$$E_M^\Delta(x,T,f):=\frac{1}M\sum_{j=1}^Mf((X_T^{\Delta,0,x})^j), $$where the $ ((X_T^{\Delta,0,x})^j)_{j\in \leftB 1,M\rightB}$ are independent copies of the scheme \eqref{EULER} starting at $x$ at time $0$ and evaluated at time $T$.
\end{trivlist}

\noindent The global error between $\E_x[f(X_T)] $, the quantity to estimate, and its implementable approximation $E_M^\Delta(x,T,f) $ can be decomposed as follows:
\begin{eqnarray}
\label{ERROR}
{\cal E}(\Delta,M,x,T,f)&:=&(\E_x[f(X_T)]-\E_x[f(X_T^\Delta)])+(\E_x[f(X_T^\Delta)]-E_M^\Delta(x,T,f))\nonumber\\
&:=&{\cal E}_D(\Delta,x,T,f)+{\cal E}_S(\Delta,M,x,T,f).
\end{eqnarray}

The term ${\cal E}_D(\Delta,x,T,f) $ corresponds to the \textit{discretization error} and has been widely investigated in the literature since the seminal work of Talay and Tubaro \cite{tala:tuba:90}. %, see also Konakov and Mammen \cite{kona:mamm:02} under the current assumption \A{A}. 
For the \textit{standard} Euler scheme, this contribution usually yields an error of order $\Delta $, provided the coefficients $b,\sigma $ and the function $f$ are sufficiently smooth, which are the assumptions required in \cite{tala:tuba:90}, or that $b,\sigma$ satisfy some non-degeneracy assumptions which allow to weaken the smoothness assumptions on $f$. 
This is for instance the case in Bally and Talay \cite{ball:tala:96:1}  who obtain the expected order for a bounded measurable $f$ and smooth coefficients satisfying a (possibly weak) hypoellipticity condition. Their proof relies on Malliavin calculus.  When the diffusion coefficient is uniformly elliptic and bounded, if $b,\sigma $ are also assumed to be three times continuously differentiable, the control at order $\Delta $ for $ {\cal E}_D(\Delta,x,T,f)$  can be derived from Konakov and Mammen \cite{kona:mamm:02} who use a more direct \textit{parametrix} approach.
When the Gaussian increments of the standard Euler scheme are replaced by more general (possibly discrete) random variables $(Y_i)_{i\ge 1}$ having the same covariance matrix and odd moments up to order 5 as  the standard Gaussian vector of $\R^q$, it can be checked that the error expansion at order $\Delta$ of \cite{tala:tuba:90} still holds for $b,\sigma,f$ smooth enough. In that framework we also mention the works of Konakov and Mammen \cite{kona:mamm:00}, \cite{kona:mamm:01}, concerning local limit theorems for the difference between \eqref{EqDfSt} and the scheme \eqref{EULER}. As in \cite{kona:mamm:02}, the coefficients are supposed to be smooth and $\sigma $ uniformly elliptic.  The associated error is then of order $\Delta^{1/2} $, speed of the Gaussian local limit theorem, see Bhattacharya and Rao \cite{bhat:rao:76}.

The term ${\cal E}_S(\Delta,M,x,T,f) $ in \eqref{ERROR} corresponds to the \textit{statistical error}. Under some usual integrability conditions, i.e. $f(X_T^\Delta)\in L^2(\P) $, it is \textit{asymptotically} controlled by the central limit theorem. A first non-asymptotic result is given by the Berry-Essen theorem provided $f(X_T^\Delta)\in L^3(\P) $, but for practical purposes, the crucial quantity to control non-asymptotically is the deviation between the empirical mean $E_M^\Delta(x,T,f) $ and the real one $\E_x[f(X_T^\Delta)]$. Precisely, for a \textit{fixed} $M$ and a given threshold $r>0$, one would like to give bounds on the quantity $\P[|E_M^\Delta(x,T,f)-\E_x[f(X_T^\Delta)]|>r] $. 

%%%%%% Attention RQ Malrieu et Talay, en fait a sigma constant on a concentration Gaussienne, le mentionner?

%For a bounded uniformly elliptic diffusion coefficient $\sigma$, assuming also that $b,\sigma $ are uniformly $\eta $-H\"older continuous in space ($\eta\in (0,1]$),
In the ergodic framework and for a constant diffusion coefficient %a Gaussian concentration property has 
Gaussian controls have been obtained by 
 Malrieu and Talay \cite{malr:tala:06}. %established a logarithmic Sobolev inequality for the scheme, thus implying \eqref{CONC_GAUSS}. 
 In the current context and for the \textit{standard} Euler scheme, a first attempt to establish two-sided Gaussian bounds for ${\cal E}_S(\Delta,M,x,T,f)  $ can be found in \cite{lema:meno:10} under some non-degeneracy conditions up to a systematic bias independent of $M$. 

In the current work we assume that the coefficients satisfy the mild smoothness condition:
 \begin{trivlist}
\item[\A{A}] The coefficients $b,\sigma$ are uniformly Lipschitz continuous in space uniformly in time, $\sigma $ is bounded.
 \end{trivlist}
Note that we do not assume any non-degeneracy condition on $\sigma$ in \A{A}.% no specific non-degeneracy on $\sigma $ is needed.
 
 We next show that when the innovations satisfy \eqref{CONC_GAUSS}, the Gaussian concentration property transfers to the statistical error $E_M^\Delta(x,T,f)-\E_x[f(X_T^\Delta)] $. In particular we get rid off the systematic bias in \cite{lema:meno:10}. 
 The key tool consists in writing the deviation using the same kind of decompositions that are exploited in \cite{tala:tuba:90} for the analysis of the discretization error. Denote by $X_{T}^{\Delta,t_i,x} $ the value at time $T$ of the scheme \eqref{EULER} starting from $x\in \R^d$ at time $t_i,\ i\in \leftB 0,N\rightB$ and by $\F_{i}:=\sigma(Y_j,\ j\le i) $ the filtration generated by the innovations. We write
 \begin{eqnarray*}
f(X_T^{\Delta,0,x})-\E[f(X_T^{\Delta,0,x})]&:=&\bsum{i=1}^{N}\E[f(X_T^{\Delta, 0,x})|\F_i]-\E[f(X_T^{\Delta, 0,x})|\F_{i-1}]\\
&=& \bsum{i=1}^N \E[f(X_T^{\Delta,0,x})|X_{t_i}^{\Delta,0,x}]-\E[f(X_T^{\Delta,0,x})|X_{t_{i-1}}^{\Delta,0,x}],
 \end{eqnarray*}
 using the Markov property for the last equality. Introducing the function $v^\Delta(t_i,x):=\E[f(X_T^\Delta)|X_{t_i}^\Delta=x], \ (i,x)\in \leftB 0,N\rightB\times \R^d $, we obtain:
 \begin{eqnarray}
f(X_T^{\Delta,0,x})-\E[f(X_T^{\Delta,0,x})]&:=&\bsum{i=1}^N v^\Delta(t_i,X_{t_i}^{\Delta,0,x})-v^\Delta(t_{i-1},X_{t_{i-1}}^{\Delta,0,x}).\label{DEC_1}
\end{eqnarray}
The definition of $v^\Delta$ now yields:
\begin{eqnarray}
f(X_T^{\Delta,0,x})-\E[f(X_T^{\Delta,0,x})]
     &=&\bsum{i=1}^N v^\Delta(t_i,X_{t_i}^{\Delta,0,x})-\E[v^\Delta(t_i,X_{t_{i}}^{\Delta,0,x})|X_{t_{i-1}}^{\Delta,0,x} ]\nonumber \\
    &=&\bsum{i=1}^N f_i^\Delta(X_{t_{i-1}}^{\Delta,0,x},\sqrt{\Delta}Y_i)-\E[f_i^\Delta(X_{t_{i-1}}^{\Delta,0,x},\sqrt{\Delta} Y_i)|X_{t_{i-1}}^{\Delta,0,x}],\nonumber\\
    \label{DEC_2}
 \end{eqnarray}
where $f_i^\Delta(x,y):=\E[f(X_T^\Delta)|X_{t_i}^\Delta=x+b(t_{i-1},x)\Delta+\sigma(t_{i-1},x) y]$, for all $(i,x,y)\in \leftB1,N\rightB\times \R^d\times \R^q$.

 The decomposition \eqref{DEC_1} is similar to the first step of the analysis of the discretization error. In that framework, $v^\Delta(t_i,X_{t_i}^{\Delta,0,x})$ is replaced
 by $v(t_i,X_{t_i}^{\Delta,0,x})=\E[f(X_T^{t_i,X_{t_i}^{\Delta,0,x}})] $, that is the expectation involving the diffusion at time $T$ starting from the current value of the scheme at $t_i$, see \cite{tala:tuba:90}. Under some non degeneracy assumptions or smoothness of the coefficients, $v$ is smooth and It\^o-Taylor expansions lead to the previously mentioned first order error for ${\cal E}_D(\Delta,x,T,f) $.
 
 To analyze the statistical error, the key idea is to exploit recursively from \eqref{DEC_2} that the increments of the scheme \eqref{EULER} satisfy \eqref{CONC_GAUSS}. The Gaussian concentration property will readily follow provided the $f_i$ are Lipschitz in the variable $y$. Under \A{A}, this smoothness is actually derived  from direct stability arguments using flow techniques, see Proposition \ref{PROP_LIP_EUL} and its proof in Section \ref{CTR_LIP}. 
 
%%%%% La il faut parler de Blower et Bolley

Let us here mention the work of Blower and Bolley \cite{blow:boll:06} who obtained Gaussian concentration properties for the joint law of the first $n$ positions of stochastic processes (possibly non Markov) with values in general separable metric spaces. This result is in some sense much stronger than ours, since it can for instance yield to non asymptotic controls of the Monte-Carlo error for smooth functionals of the path, such as the maximum. However, some continuity assumptions in Wasserstein metric are assumed on the transition measures of the process, see e.g. condition \textit{(ii)} in their Theorems 1.2 and 2.1. This is required from the coupling techniques used in the proof. Checking this kind of continuity can be hard in practice, in \cite{blow:boll:06} the authors give some sufficient conditions that require the transition laws to be absolutely continuous and smooth, see their Proposition 2.2.  In the current work we only need the property \eqref{CONC_GAUSS} for the innovations, which can in particular hold for discrete laws.

%Also, the global control on the joint law in \cite{blow:boll:06} ignores the smoothing effects of non degenerated heat kernels that can be useful to establish a Gaussian concentration for non smooth terminal functions $f$ such as the indicator of a set, which is of importance in finance for digital options, or fatigue analysis of materials.  

Also, we want to stress that, even if the concentration results coincide when the innovations $(Y_i)_{i\in \N^*}$ have a smooth density, the nature of the proofs is different. Blower and Bolley exploit optimal transportation techniques whereas our approach consists in adapting the PDE arguments used for the analysis of the discretization error to the current setting. It is actually striking that a similar error decomposition can be used for investigating both the discretization and statistical error. %Furthermore, the same kind of assumptions appear, i.e. smoothness of the coefficients or non degeneracy conditions, in order to guarantee the smoothness of the functions involved.

%We conclude mentioning as well the results of Bolley, Guillin, Villani \cite{boll:guil:vill:07} and Boissard \cite{bois:11}, who obtained non-asymptotic deviation bounds in the 1-Wassertstein metric between a reference measure $\nu$ and its empirical version $L_{n}:=\frac{1}{n}\sum_{i=1}^{n} \delta_{X_i}$ where $(X_{i})_{ i\in \leftB 1,  n\rightB}$ is an i.i.d. sequence of random variables with common law $\nu$.  These results respectively  rely  on a non-asymptotic version of Sanov's Theorem and concentration inequalities. Though these results are of a different nature, they imply somehow stronger concentration inequalities such that 
%$$
%\P\left(\sup_{f, 1-Lip} \left|\frac{1}{M} \sum_{j=1}^{M} f(X_j)- \E[f(X)] \right| \geq r\right) \leq C(r)\exp\left(-K M r^{2}\right)
%$$
%\noindent where the constants $C(r)$ and $K$ may be explicitely computed. This kind of uniform deviation bound are of interest in statistics and numerical probability from a practical point of view. We do not intend to develop these aspects but similar bounds could be established in our context.

We conclude mentioning some works related to the deviations of the 1-Wasserstein distance between a reference measure and its empirical version.
In the i.i.d. case,  such results were first obtained for different concentration regimes by Bolley, Guillin, Villani \cite{boll:guil:vill:07} relying on a non-asymptotic version of Sanov's Theorem. Some of these results have also been derived by Boissard \cite{bois:11} using concentration inequalities and extended to ergodic Markov chains
 up to some contractivity assumptions in the Wasserstein metric on the transition kernel. In the i.i.d. case  and Gaussian concentration regime, these results lead to the following type of estimates:
$$
\P\left(\sup_{f, 1-Lip} \left|\frac{1}{M} \sum_{j=1}^{M} f(Z_j)- \E[f(Z)] \right| \geq r\right) \leq C(r)\exp\left(-K M r^{2}\right)
$$
\noindent where the $(Z_j)_{j\in \N^*} $ are i.i.d. having the same distribution as $Z$  and the constants $C(r)$ and $K$ may be explicitly but tediously computed. This kind of uniform deviation bounds are of interest in statistics and numerical probability from a practical point of view.
They can indeed lead to deviation bounds for the estimation of the density of the invariant measure of a Markov chain, see \cite{boll:guil:vill:07}. 
However, the (possibly large) constant $C(r)$ is the trade-off to obtain uniform deviations over all Lipschitz functions.
 We do not intend to develop these aspects but similar bounds could be established in our context. 

 \subsection{Robbins-Monro Stochastic Approximation Algorithm}
 \label{RM_SEC}
Besides our considerations for the Euler scheme, we derive non asymptotic bounds for stochastic approximation algorithms of Robbins-Monro type. %They stand for a family of simulation based 
These recursive %optimization 
algorithms aim at finding a zero of a continuous function $h:\R^{d} \rightarrow \R^{d}$ which cannot be directly computed but only estimated through simulation. Such procedures are commonly used in a convex optimization framework since minimizing a function amounts to finding a zero of its gradient. %Hence, this refers to estimating 
Precisely, the goal is to find a solution $\theta^{*}$ to $h(\theta):=\E[H(\theta,Y)]=0$, where $H: \R^{d} \times \R^{q} \rightarrow \R^{d}$ is a Borel function and $Y$ is a given $\R^{q}$-valued random variable. 
 Even though $h(\theta)$ cannot be directly computed,  it is assumed that the random variable $Y$ can be easily simulated (at least at a reasonable cost), and  also that $H(\theta,y)$ can be easily computed for any couple $(\theta,y) \in \R^{d}\times \R^{q}$. 
%It is assumed that though we cannot directly compute $h(\theta)$, we can easily simulate the distribution $\mu$ of $Y$, at least at a reasonable cost, and we can compute $H(\theta,y)$ for any couple $(\theta,y) \in \mathbb{R}^{d}\times \mathbb{R}^{q}$. 
%Then, one defines the following recursive scheme
The Robbins-Monro algorithm is the following recursive scheme
\begin{equation}
\label{RM}
\theta_{n+1} = \theta_{n} - \gamma_{n+1} H(\theta_{n},Y_{n+1}), \ n \geq 0, \ \theta_{0} \in \R^{d},
\end{equation}
%%% j'ai enleve with law \mu dans la mesure ou l'on n'y fait plus reference.
\noindent where $(Y_{n})_{n\geq1}$ is an i.i.d. $\R^{q}$-valued sequence of  random variables %with law $\mu$ 
defined on a probability space $(\Omega, \mathcal{F}, \mathbf{P})$ and $(\gamma_{n})_{n\geq1}$ is a sequence of non-negative deterministic steps satisfying the usual assumption
\begin{equation}
\label{STEP}
\sum_{n\geq1} \gamma_{n} = + \infty, \ \ \mbox{and} \ \ \sum_{n\geq1} \gamma_{n}^{2} < + \infty.
\end{equation} 
%%%% De S pour N: Phrase de mon cru... Je crois pas dire trop de betises...
\noindent When the function $h$ is the gradient of a potential, the iterative scheme \eqref{RM} can be viewed as a stochastic gradient algorithm. Indeed, replacing $H(\theta_n,Y_{n+1})$ by $h(\theta_n)$ in \eqref{RM} leads to the usual deterministic gradient method. One of the ideas in \eqref{RM} is to take advantage of an averaging effect along the scheme due to the specific form of $h(\theta):=\E[H(\theta,Y)]$. This allows to avoid the explicit computation or estimation of $h$.
\noindent We refer to \cite{Duflo1996}, \cite{Kushner2003} for some general convergence results of the sequence $(\theta_{n})_{n\geq0}$ defined by \eqref{RM} towards its target $\theta^{*}$ under the existence of a so-called \emph{Lyapunov function}, $i.e.$ a continuously differentiable function $L:\mathbf{R}^{d}\rightarrow \mathbf{R}_{+}$ such that $\nabla L$ is Lipschitz, $|\nabla L|^{2} \leq C(1+L)$ for some positive constant $C$ and 
$$
 \left\langle \nabla L , h \right\rangle \geq 0.
$$

\noindent See also \cite{Laruelle2012} for a convergence theorem under the existence of a \emph{pathwise Lyapunov function}. In the sequel, it is assumed that $\theta^{*}$ is the unique solution of the equation $h(\theta)=0$ and that $(\theta_{n})$ defined by \eqref{RM} converges $a.s.$ towards $\theta^{*}$.

We assume that the innovations $(Y_i)_{i\in \N^*} $ satisfy \eqref{CONC_GAUSS} for some $\alpha>0$ and also that the following conditions on the function $H$ and the step sequence $(\gamma_{n})_{n\geq1}$ in \eqref{RM} are in force:
\begin{trivlist}
\item[\A{HL}]The map $(\theta,y)\in \R^{d} \times \R^{q} \mapsto H(\theta,y)$ is uniformly Lipschitz continuous.
\item[\A{HUA}]The map $h:\theta \in \R^d \mapsto \E[H(\theta,Y)]$ is continuously differentiable in $\theta $ and there exists $\underline{\lambda}>0 $ s.t. $\forall \theta \in \R^{d}, \ \forall \xi\in \R^d, \  \underline{\lambda} |\xi|^2 \le \langle Dh(\theta)\xi,\xi\rangle %\le \overline{\lambda} |\xi|^2 
$ 
(\textit{Uniform Attractivity}). 
\end{trivlist}
 In order to derive a Central Limit Theorem for the sequence $(\theta_{n})_{n\geq1}$ as described in \cite{Duflo1996} or \cite{Kushner2003}, it is commonly assumed that the matrix $Dh(\theta^{*})$ is \emph{uniformly attractive}. In our current framework, this local condition on the Jacobian matrix of $h$ at the equilibrium is replaced by the uniform assumption \A{HUA}. %which seems quite natural since we want to 
This allows to derive non-asymptotic concentration bounds %for $|\theta_{n}-\theta^*|-\E[|\theta_{n}-\theta^*|], \ n\ge 1$ 
uniformly w.r.t. the starting point $\theta_0$. %The sensitivity of the scheme \eqref{RM} w.r.t. to the initial point is only reflected in the bias $\delta_n:=\E[|\theta_{n}-\theta^*|] $.

Note that under \A{HUA} and the \emph{linear growth assumption}
$$
\forall \theta \in \R^{d}, \ \ \E\left[\left|H(\theta,Y)\right|^{2}\right] \leq C(1+|\theta-\theta^{*}|^{2}),
$$

\noindent (which is satisfied if \A{HL} holds and $Y \in L^{2}(\mathbf{P})$) the function $L:\theta \mapsto \frac{1}{2}\left|\theta - \theta^{*}\right|^{2}$ is a Lyapunov function for the recursive procedure defined by \eqref{RM} so that one easily deduces that $\theta_{n} \rightarrow \theta^{*}$, $a.s.$ as $n\rightarrow + \infty$.

%%%%% Commenter un peu plus (UHA) comme indique.

As for the Euler scheme, we decompose the global error between the stochastic approximation procedure $\theta_{n}$ at a given time step $n$ and its target $\theta^{*}$ as follows:
\begin{eqnarray}
z_{n}  := \left|\theta_{n} - \theta^{*}\right| & = & \left(\left|\theta_{n} - \theta^{*}\right| - \E[\left|\theta_{n} - \theta^{*}\right|]  \right) + \E[\left|\theta_{n} - \theta^{*}\right|] \nonumber\\
      &  := & \mathcal{E}_{Emp}(\gamma,n,H,\underline{\lambda}, \alpha) + \delta_{n}\label{DEC_RM}
\end{eqnarray}

\noindent where $\delta_n := \E[\left|\theta_{n} - \theta^{*}\right|]$.

The term $\mathcal{E}_{Emp}(\gamma,n,H,\underline{\lambda}, \alpha)$ corresponding to the difference between the absolute value of the error at time $n$ and its mean can be viewed as an \emph{empirical error}. As for the Euler scheme, the Gaussian concentration property transfers to this quantity under \A{HL} and \A{HUA}. The strategy consists in introducing again a telescopic sum of conditional expectations. Denoting for all $i\in \N,\ \F_i:=\sigma(Y_j,\ j\le i)$ (i.e. $(\F_i)_{i\in \N}$ is the natural filtration of the algorithm), we write for all $n\in \N^*$:
\begin{eqnarray*}
\mathcal{E}_{Emp}(\gamma,n,H,\underline{\lambda}, \alpha)= |z_n| - \E[|z_n|] & = & \bsum{i=1}^n \E[|z_n||\F_{i}]-\E[|z_n||\F_{i-1}]  \\
&= &\bsum{i=1}^{n } v_i(\theta_i)-\E[v_i(\theta_i)|\F_{i-1}], \\
       & = & \bsum{i=1}^n f_i^\gamma(\theta_{i-1},Y_i)-\E[f_i^\gamma(\theta_{i-1},Y_i)|\F_{i-1}],
\end{eqnarray*}

\noindent where we used the Markov property for the second equality and we introduced the notations $v_i(\theta):=\E[|\theta_n-\theta^*||\theta_i=\theta],\ \forall (i,\theta)\in \leftB 1,n \rightB\times \R^d,\ f_i^\gamma(\theta,y)=v_i(\theta -\gamma_i H(\theta,y))$. The stability of the Gaussian concentration property is then derived using that the $f_i^\gamma$ are Lipschitz in the variable $y$, see Proposition \ref{CTR_FIG}. 

The term $\delta_n$ in \eqref{DEC_RM} corresponds to the \emph{bias} of the sequence $(\theta_{n})_{n\geq0}$ with respect to its target $\theta^{*}$. This contribution strongly depends on the choice of the step sequence $(\gamma_{n})_{n\geq1}$ and the initial point $\theta_0 $. Under \A{HL} and \A{HUA}, we analyze this quantity in Proposition \ref{CTR_BIAS}.

%%%%%% Version S
%Set now for all $n\in \N, z_n:=\theta_n-\theta^* $. 
%As in the previous section, we want to control the deviations of $|z_n|$. We use the same strategy introducing a telescopic sum of conditional expectations. For all $n\in \N^* $, write
%\begin{eqnarray*}
%|z_n| - \E[|z_n|] & = & \bsum{i=1}^n \E[|z_n||\F_{i}]-\E[|z_n||\F_{i-1}]  = \bsum{i=1}^{n } v_i(\theta_i)-\E[v_i(\theta_i)|\F_{i-1}], \\
%       & = & \bsum{i=1}^n f_i^\gamma(\theta_{i-1},Y_i)-\E[f_i^\gamma(\theta_{i-1},Y_i)|\F_{i-1}],
%\end{eqnarray*}

%\noindent where we used the Markov property for the second equality and we introduced the notations $v_i(\theta):=\E[|\theta_n-\theta^*||\theta_i=\theta],\ \forall (i,\theta)\in \leftB 1,n \rightB\times \R^d,\ f_i^\gamma(\theta,y)=v_i(\theta -\gamma_i H(\theta,y)) $. The previous assumptions \A{HL}, \A{HUA} allow to control the Lipschitz constants in $y $ of the $f_i^\gamma$, see Proposition \ref{CTR_FIG}.

 \mysection{Main Results}
 
 \subsection{Deviations on the Euler Scheme}
 \begin{THM}[Concentration Bounds for the Euler scheme]
\label{THM_EULER}
   Denote by $X_T^\Delta$ the value at time $T$ of the scheme \eqref{EULER} associated to the diffusion \eqref{EqDfSt}. 
Assume that  the innovations $(Y_i)_{i\in \N^*} $ in \eqref{EULER} satisfy \eqref{CONC_GAUSS} for some $\alpha>0 $ and that 
the coefficients $b,\sigma$ satisfy \A{A}.  Let $f$ be a real valued uniformly Lipschitz continuous function on $\R^d$. 
  For all $M\in \N^*$ and all $r\ge 0$, one has
\begin{eqnarray*} 
 &&\P_x[|\frac{1}{M}\bsum{i=1}^M f((X_T^\Delta)^i)-\E_x[f(X_T^\Delta)]|\ge r]\le 2\exp(-\frac{r^2M}{T\Psi(T,f,b,\sigma,q)}),\\
 && \Psi(T,f,b,\sigma,q):=4\alpha[f]_1^2|\sigma|_\infty^2\exp\left(2([b]_1+c[\sigma]_1(1\vee c[\sigma]_1) ) T\right),
\end{eqnarray*}
\noindent where $q$ is the dimension of the underlying Brownian motion in \eqref{EqDfSt} and $c:=c(q)$.
 \end{THM}

Note that in the above theorem, we do not need any non-degeneracy condition on the diffusion coefficient. As developed in Section \ref{DESCR_EUL}, see \eqref{DEC_2}, to handle the previous quantity we rewrite $f(X_T^{\Delta})-\E[f(X_T^{\Delta})]:=\sum_{i=1}^N f_i^\Delta(X_{t_{i-1}}^\Delta,\sqrt \Delta Y_{i})-\E[f_i^\Delta(X_{t_{i-1}}^\Delta,\sqrt \Delta Y_{i})|\F_{i-1}]$, where  $f_i^\Delta(x,y):=\E[f(X_T^\Delta)|X_{t_i}^\Delta=x+b(t_{i-1},x)\Delta+\sigma(t_{i-1},x) y]$, for all $(i,x,y)\in \leftB1,N\rightB\times \R^d\times \R^q$.
If at some point along the time-discretization the process has a degenerate diffusion term, we can see that the difference $f_i^\Delta(X_{t_{i-1}}^\Delta,\sqrt \Delta Y_{i})-\E[f_i^\Delta(X_{t_{i-1}}^\Delta,\sqrt \Delta Y_{i})|\F_{i-1}] $
will not give any additional contribution in the global deviation.
%since we are considering deviations from the mean

With respect to the previous work \cite{lema:meno:10}, we got rid off the systematic bias. Anyhow, the concentration constants now depend on the Lipschitz constant of the function $v^\Delta(0,x):=\E[f(X_T^\Delta)|X_0^\Delta=x]$ which has order $\Psi(T,f,b,\sigma,q)^{1/2}$. This magnitude corresponds to the product  of the  Lipschitz constant of the final function $f$ and the mean of the Lipschitz constant for the flow of the scheme, which gives the exponential dependence in time, see Proposition \ref{PROP_LIP_EUL} and its proof for details.

\begin{REM}[Extension to smooth functionals of the path] We point out that the previous concentration results could be extended to some smooth functionals of the path such as the maximum for a scalar scheme. Indeed, introducing in that case the additional state variable $(M_{t_i}^\Delta)_{i\in \N}:=(\max_{j\in \leftB 0,i\rightB}X_{t_j}^\Delta)_{i\in \N} $, the couple $(X_{t_i}^\Delta,M_{t_i}^\Delta)_{i\in \N} $ is Markovian and the flow arguments of Proposition \ref{PROP_LIP_EUL} could be extended to the couple for Lipschitz functions in both variables.
\end{REM}

\begin{REM}[Linear SDEs and concentration]
Observe that it is the boundedness of $\sigma$ that gives the Gaussian concentration regime. However, in many popular models in finance, the diffusion coefficient is linear, see e.g. the Black-Scholes like dynamics $X_t=x_0+\int_0^tb(X_s)X_sds +\int_0^t\sigma(X_s)X_sdW_s $ for smooth, bounded coefficients $b,\sigma$. For the estimation of $\E[f(X_T^\Delta)] $ of the associated Euler scheme, if $f$ is bounded then the Gaussian concentration holds for the statistical error from the Bolley and Villani criterion applied to $f(X_T^\Delta)$. However, for a general Lipschitz function, the expected concentration is the log-normal one. 
\end{REM}

%
%
%
%  Faire la remaruqe sur le cas lineaire, comme suggere par le referee 2.
\subsection{Deviations for Robbins-Monro algorithms}
 \begin{THM}[Concentration Bounds for Robbins-Monro algorithms]
 \label{THM_RM}
 Assume that the function $H$ of the recursive procedure $(\theta_{n})_{n\geq0}$ (with starting point $\theta_{0} \in \R^{d}$) defined by \eqref{RM} satisfies \A{HL} and \A{HUA}, and that the step sequence $(\gamma_{n})_{n\geq1}$ satisfies \eqref{STEP}. Suppose that the law of the innovation satisfies \eqref{CONC_GAUSS}, $\alpha>0 $.% and $\int |x|^{2} d\mu(x) < + \infty$. Fix $N \in \N^{*}$. %%%%% Ben come on a concentration gaussienne, la variance exist bien!
 Then, for all $N\in \N^*$ and all $r \geq 0$,
$$
\P\left(\left|\theta_{N} - \theta^{*}\right| \geq r + \delta_{N}\right) \leq \exp\left(- \frac{r^{2}}{\alpha [H]^{2}_{1}\Pi_{N}\sum_{k=1}^{N}\gamma^{2}_{k}/\Pi_{k}}\right)
$$

\noindent where $\Pi_{N}:=\prod_{k=0}^{N-1} \left( 1 - 2 \underline{\lambda} \gamma_{k+1} + [H]^{2}_{1} \gamma^{2}_{k+1}\right)$ and $\delta_{N}:=\E\left[\left|\theta_{N}-\theta^{*}\right|\right]$. Moreover, the bias $\delta_{N}$ at step $N$ satisfies
\begin{eqnarray*}
\delta_{N} \leq \exp\left(-\underline{\lambda} \Gamma_{N}\right) \left|\theta_0 - \theta^{*}\right| +  [H]_{1} \sigma_{Y} \left( \sum_{k=0}^{N-1} e^{-2\underline{\lambda} (\Gamma_{N} - \Gamma_{k+1})} \gamma^2_{k+1}\right)^{\frac{1}{2}},
\end{eqnarray*}

\noindent where $\Gamma_{N} :=\sum_{k=1}^{N} \gamma_{k}$, $\sigma_{Y} := \E\left[F^{2}(Y)\right]^{1/2}<+\infty$, with $F:y \mapsto \E\left[|y-Y|\right]$. % and $c_{1}:=c(q)>1$ ($q$ being the dimension of $Y$) is a constant.
\end{THM}

Concerning the choice of the step sequence $(\gamma_{n})_{n\geq1}$ and its impact on the concentration rate and bias, we obtain the following results:
\begin{itemize} 
\item If we choose $\gamma_{n}=\frac{c}{n}$, with $c>0$. Then $\delta_{N} \rightarrow 0$, $N\rightarrow+\infty$, $\Gamma_{N}=c \log(N) + c'_{1} + r_{N}$, $c'_{1}>0$ and $r_{N}\rightarrow 0$, so that $\Pi_{N} =\O( N^{-2c \underline{\lambda}})$.
\begin{itemize}

\item If $c < \frac{1}{2 \underline{\lambda}}$, the series $\sum_{k=1}^{N}\gamma^{2}_{k}/\Pi_{k}$ converges so that we obtain 
$\Pi_{N}\sum_{k=1}^{N}\gamma^{2}_{k}/\Pi_{k} = \mathcal{O}(N^{-2c \underline{\lambda}})$.

\item If $c > \frac{1}{2 \underline{\lambda}}$, a comparison between the series and the integral yields $\Pi_{N}\sum_{k=1}^{N}\gamma^{2}_{k}/\Pi_{k}$ $= \O(N^{-1})$.
\end{itemize}

\medskip

\noindent Let us notice that we find the same critical level for the constant $c$ as in the Central Limit Theorem for stochastic algorithms. Indeed, if $c> \frac{1}{2\mathcal{R}e(\lambda_{min})}$ where $\lambda_{min}$ denotes the eigenvalue of $Dh(\theta^{*})$ with the smallest real part then we know that a Central Limit Theorem holds for $(\theta_{n})_{n\geq1}$ (see e.g. \cite{Duflo1996}). However, this local condition on the Jacobian matrix of $h$ at the equilibrium is replaced by a uniform assumption in our framework. This is quite natural since we want to derive non-asymptotic bounds for the stochastic approximation \eqref{RM}.

Concerning the bias we have the following bound:
$$
\delta_N\le \frac {|\theta_0-\theta^*|} {N^{\underline{\lambda} c}}+ K [H]_{1} \sigma_{Y} \frac{1}{N^{\underline{\lambda}c \wedge \frac{1}{2}}},\ K:=K(c). 
$$

\item If we choose $\gamma_{n}=\frac{c}{n^{\rho}}$, $c>0$, $\frac{1}{2} < \rho < 1$, then $\delta_{N} \rightarrow 0$, $\Gamma_{N} \sim \frac{c}{1-\rho} N^{1-\rho}$ as $N \rightarrow + \infty$ and elementary computations show that there exists $C>0$ s.t. for all $N\ge 1$, $\Pi_N\le C\exp(-2\underline{\lambda}\frac{c}{1-\rho}N^{1-\rho}) $. Hence, for all $\epsilon\in (0,1-\rho) $ we have:
\begin{eqnarray*}
\Pi_{N}\sum_{k=1}^{N} \gamma_k^2 \Pi^{-1}_{k}&\le & c^2 \left\{\Pi_{N}\Pi_{N-N^{\rho+\epsilon}}^{-1}\sum_{k=1}^{N-N^{\rho+\epsilon}} \frac1{k^{2\rho}}+\sum_{k=N-N^{\rho+\epsilon}+1}^{N} \frac1{k^{2\rho}}\right\}\\&\le &  c^2 \left\{C\exp(-2\underline{\lambda}\frac{c}{1-\rho}(N^{1-\rho}-(N-N^{\rho+\epsilon} )^{1-\rho})) \right.\\
&&\left.+  \frac{N^{\rho+\epsilon}}{(N-N^{\rho+\epsilon}+1)^{2\rho}}\right\}\\
&\le & c^2 \left\{C\exp(-2\underline{\lambda}c N^{\epsilon})+  \frac1{N^{\rho-\epsilon}} \right\}.
\end{eqnarray*}

\noindent Up to a modification of $\epsilon$, this yields $\Pi_{N}\sum_{k=1}^{N} \gamma_k^2 \Pi^{-1}_{k}=o(N^{-\rho+\epsilon}),\ \epsilon\in (0,1-\rho) $. 

\end{itemize}

Concerning the bias, from the above control, we have the following bound:
\begin{eqnarray*}
\delta_N\le \exp\left( -\frac{\underline{\lambda}c }{1-\rho}  N^{1-\rho} \right)|\theta_0-\theta^*|+  [H]_{1} \sigma_{Y} \frac{K}{N^{\frac{\rho }{2}- \epsilon }},%\\
 K:=K(c),\ \forall \epsilon>0. 
\end{eqnarray*}

%The main difference with the previous choice of time step is that, 
Since each step is bigger compared to the case $\gamma_{n}=\frac{c}{n}$, the impact of the initial difference $|\theta_0-\theta^*| $ is exponentially small. 
%On the other hand, the other contributions are bigger.
                            
%
%
%
%
 \mysection{Abstract concentration properties for a general evolution scheme}
In this section we assume that $(Y_i)_{i\in \N^*}$ is a sequence of i.i.d. $\R^q$-valued random variables whose law $\mu$ satisfies the Gaussian concentration property  \eqref{CONC_GAUSS} for a given $\alpha>0 $. 
\begin{PROP}[Gaussian concentration for a stochastic evolution scheme]
\label{HERBST_COND}
Fix $N\in \N^*$. Define for all $i \in \leftB 1, N \rightB$, ${\scD}_{i}:=f_i(X_{i-1},Y_{i})-\E[f_i(X_{i-1},Y_{i})|\F_{i-1}]$ for some $\F_{i-1}$-measurable random variables $X_{i-1}$ where the real-valued functions $(f_i)_{i\in \leftB 1,N\rightB}$ are Lipschitz continuous in the $y$ variable with constants $([f_i]_{1})_{i\in \leftB 1,N\rightB}>0$ uniformly in $x$. Let $(\gamma_{i})_{i\in \leftB 1,N\rightB} $ be a given sequence of time steps. For all $r\ge 0$, we have:
$$
\P[\sum_{i=1}^{N}\gamma_i \scD_{i} \ge r]\le \exp\left(-\frac{r ^2}{\alpha \sum_{i=1}^{N}([f_{i}]_{1} \gamma_{i})^2}\right). 
$$
\end{PROP}

\textit{Proof.} 
%As in the proof of the Gaussian concentration for a measure satisfying a logarithmic Sobolev inequality, see e.g Ledoux \cite{ledo:99}, 
Set ${\cal P}(r):=\P[\sum_{i=1}^{N}\gamma_i \scD_{i} \ge r]$.
For $\lambda \geq 0$ to be specified later on, the Tchebychev exponential inequality yields:
\begin{eqnarray}
{\cal P}(r)&\le &\exp(-\lambda r)\E[\exp\left(\lambda \left[ \sum_{i=1}^{N} \gamma_{i} \scD_{i}\right] \right)]\nonumber \\
&\le & \exp(-\lambda r)\E[\exp\left(\lambda \sum_{i=1}^{N-1} \gamma_{i} \scD_{i}\right) \E\left[ \exp(\lambda\gamma_N \scD_N  ) |\F_{N-1}\right] ].
\label{COND_EC}
\end{eqnarray}

\noindent Observe now that working with regular conditional expectations, we have 
\begin{eqnarray*}
\E\left[ \exp(\lambda \gamma_N \scD_N) |\F_{N-1}\right]\nonumber\\ 
=\left.\E_\mu\left[ \exp \left(\lambda \gamma_N   (f_{N}(x,Y)-\E_\mu[f_{N}(x,Y)])  \right)\right]\right|_{ x=X_{N-1}},%\nonumber\\
%\le \left.\E_\mu\left[ \exp \left(2 \lambda b (f_{N-1}(y,X)-\E_\mu[f_{N-1}(y,X)]) \right)\right]\right|_{b=b_{N-1},c=c_{N-1}, y=Y_{N-1}}^{1/2}\nonumber \\
%\times\left.\E_\mu\left[ \exp \left(2 \lambda c (g_{N-1}(y,X)-\E_\mu[g_{N-1}(y,X)]) \right)\right]\right|_{b=b_{N-1},c=c_{N-1}, y=Y_{N-1}}^{1/2},
\end{eqnarray*}

\noindent where $Y$ is a random variable with law $\mu$. From \eqref{CONC_GAUSS}, we derive
\begin{eqnarray*}
\E\left[ \exp(\lambda \gamma_N\scD_N) |\F_{N-1}\right]
\le 
\exp(\alpha ([f_{N}]_{1} \gamma_N\lambda)^2/4 ) .
\end{eqnarray*}
 
\noindent Plugging this estimate in \eqref{COND_EC} and iterating the procedure we derive
$$
{\cal P}(r) \le \exp(-\lambda r)\exp\left(\frac{\alpha \lambda^2}{4} \sum_{i=1}^N ([f_i]_1 \gamma_i)^2\right),
$$
 
\noindent and optimizing w.r.t $\lambda $, we obtain: ${\cal P}(r) \leq \exp\left(-\frac{r^2}{\alpha \sum_{i=1}^N([f_i]_1 \gamma_i)^2}\right)$.
\mysection{Euler Scheme: Proof of the Main Results}
\label{SEC_EUL_PR}

 In order to apply Proposition \ref{HERBST_COND} from the decomposition \eqref{DEC_2}, all we need is to have a control on the Lipschitz modulus in the variable $y$ of the functions  $f_i^\Delta(x,y) $, uniformly in $x$.
 \noindent Under the current assumptions of Theorem \ref{THM_EULER}, we have the following Proposition which is proved in Section \ref{CTR_LIP}.
 \begin{PROP}[Control of the Lipschitz constants]
 \label{PROP_LIP_EUL}
 \noindent %Assume the coefficients $b,\sigma $ in \eqref{EULER} satisfy \A{A} and D
 Denote the respective Lipschitz constants of $b,\sigma $ in \eqref{EULER} by $ [b]_1,[\sigma]_1$.   
 Denote the supremum of $\sigma $ by $|\sigma|_\infty $. %Suppose also that $f$ is uniformly Lispchitz continuous and that the $(Y_i)_{i\in \N^*} $ belong to $L^2(\P) $.
 Then for all $i\in \leftB 1,N \rightB$, the functions $f_i^\Delta  $ introduced after \eqref{DEC_2} are uniformly Lipschitz continuous in the space variable $y$ uniformly in $x$ and we have that there exists $c:=c(q)$ (dimension of the underlying Brownian motion) s.t:
$$
[f_i^\Delta]_1:=\sup_{x\in \R^d, y\neq y'} \frac{|f_i^\Delta(x,y)-f_i^\Delta(x,y')|}{|y-y'|}\le  2 [f]_1|\sigma|_\infty\exp\left(\left\{ [b]_1+c[\sigma]_1(1\vee c[\sigma]_1) \right\} (T-t_i)\right).
$$

\noindent where $[f]_1 $ stands for the Lipschitz constants of the function $f$.
\end{PROP}
 
 Set $\gamma_i=1,\ \scD_{i}=f_i^{\Delta}(X_{t_{i-1}}^\Delta,\sqrt \Delta Y_{i})-\E[f_i^{\Delta}(X_{t_{i-1}}^\Delta,\sqrt \Delta Y_{i})|\F_{t_{i-1}}], \ \forall i\in \leftB 1,N\rightB$. Since the random variable $\sqrt \Delta Y_{i}$ satisfies the Gaussian concentration property $(GC(\Delta \alpha)) $, %\eqref{CONC_GAUSS} 
 %with constant  $\Delta \alpha$ 
 we derive from Propositions \ref{HERBST_COND} and \ref{PROP_LIP_EUL}:
\begin{eqnarray*}
\label{ONE_PATH_DEV}
 \P_x[f(X_T^{\Delta})-\E_x[f(X_T^\Delta)]\ge r ]&\le& \exp(-\frac{r^2}{\Delta \alpha\sum_{i=1}^{N} [f_i^\Delta]_1^2 })\nonumber\\
 &\le &\exp(-\frac{r^2}{4\alpha T [f]_1^2|\sigma|_\infty^2\exp\left(2([b]_1+c[\sigma]_1(1\vee c[\sigma]_1) )T\right)   })\nonumber\\
 &:=&\exp(-\frac{r^2}{T \Psi(T,f,b,\sigma,q)}).
\end{eqnarray*}

\noindent 
Hence the random variable $X_T^{\Delta}$ satisfies $(GC(\beta))$ for $ \beta:=T \Psi(T,f,b,\sigma,q)/[f]_1^2 $.
The bound of Theorem \ref{THM_EULER} now follows from a simple tensorization argument for independent random variables satisfying the Gaussian concentration property.
Namely, for $\lambda, r\ge 0 $,
\begin{eqnarray*}
\P_x[\frac{1}{M}\bsum{j=1}^M f((X_T^\Delta)^j)-\E_x[f(X_T^\Delta)] \ge r]\\
\le \exp(-\lambda r)\E_x\biggl[\exp\biggl(\lambda \frac{1}M (\bsum{j=1}^M f((X_T^\Delta )^j)-\E_x[f(X_T^\Delta)^j]) \biggr)\biggr]\\
\overset{(GC(\beta))}{\le}\exp(-\lambda r)\left[\exp(\frac{\beta \lambda^2[f]_1^2}{4 M^2}) \right]^M=\exp(-\lambda r+\frac{T\Psi(T,f,b,\sigma,q)\lambda ^2}{4M}).
\end{eqnarray*}
An optimization in $\lambda $ gives the result.

% 
%To derive the deviation bound between the empirical mean and the real one, the key idea is to exploit the tensorization of the concentration property \eqref{CONC_GAUSS} for independent random variables, see e.g. Ledoux \cite{ledo:99}. %%%%%%% Detailler et renvoyer
%Precisely, considering independent sequences $\bigl((Y_i^j)_{i\ge 1 }\bigr)_{j\in \leftB 1,M\rightB} $ of $\R^q $-valued i.i.d random variables whose law satisfies \eqref{CONC_GAUSS} with constant $ \alpha$, the $\R^{Mq} $-valued random variable ${\mathbf Y}_i:=(Y_i^1,\cdots, Y_i^M),\ i\in \N^* $ satisfies \eqref{CONC_GAUSS}
%with the same constant $\alpha $.
%Set now $F_M:\x =(\x_1,\cdots,\x_M)\in \R^{Md}\mapsto F_M(\x)=\frac{1}{\sqrt M}\sum_{i=1}^M f(\x_i) $. Denoting by $|.| $ the Euclidean norm of $\R^{Md}$, the Cauchy-Schwarz inequality directly gives that $[F_M]_1:=\sup_{\x\neq \x'}\frac{\frac{1}{\sqrt M}\sum_{i=1}^{M}|f(\x_i)-f(\x_i')|}{|\x-\x'|}\le [f]_1 $.  For all $r\ge 0$, 
%\begin{eqnarray*}
%\P[ \frac{1}{M}\sum_{i=1}^M f((X_T^{\Delta,0,x})^i)-\E_x[f(X_T^{\Delta,0,x})] \ge r]%\\
% &= & \P[\frac{1}{\sqrt M} \sum_{i=1}^{M} (f((X_T^{\Delta,0,x})^i)-\E[f((X_T^{\Delta,0,x})^i)])\ge r\sqrt M ]\\
%& = & \P[F_M(\X_T^\Delta)-\E[F_M(\X_T^{\Delta,0,x})]\ge r\sqrt M],
%\end{eqnarray*}
% \noindent where $\X_T^{\Delta,0,x} :=((X_T^{\Delta,0,x})^1,\cdots,(X_T^{\Delta,0,x})^M)$. 
% The claim then follows as for $M=1$ from Proposition \ref{HERBST_COND}.
%\noindent which proves the Theorem. 
%
%
%
\mysection{Robbins-Monro Algorithm: Proof of the Main Results}
\label{PROOF_RM}
%%%% J'ai deplace la demarche dans la premiere section pour le paralele avec le schema d'Euler

%Set first for all $n\in \N, z_n:=\theta_n-\theta^* $. 
%As in the previous section, we want to control the deviations of $|z_n|$. We use the same strategy introducing a telescopic sum of conditional expectations. For all $n\in \N^* $, write
%\begin{eqnarray*}
%|z_n| - \E[|z_n|] & = & \bsum{i=1}^n \E[|z_n||\F_{i}]-\E[|z_n||\F_{i-1}]  = \bsum{i=1}^{n } v_i(\theta_i)-\E[v_i(\theta_i)|\F_{i-1}], \\
%       & = & \bsum{i=1}^n f_i^\gamma(\theta_{i-1},Y_i)-\E[f_i^\gamma(\theta_{i-1},Y_i)|\F_{i-1}],
%\end{eqnarray*}

%\noindent where we used the Markov property for the second equality and we introduced the notations $v_i(\theta):=\E[|\theta_n-\theta^*||\theta_i=\theta],\ \forall (i,\theta)\in \leftB 1,n \rightB\times \R^d,\ f_i^\gamma(\theta,y)=v_i(\theta -\gamma_i H(\theta,y)) $.

With the notations of Section \ref{RM_SEC}, in order to apply Proposition \ref{HERBST_COND} we have to control the Lipschitz constants in $y$
of the functions $f_i^\gamma(\theta,y)=v_i(\theta -\gamma_i H(\theta,y)), \forall (i,\theta,y)\in \leftB 1,n \rightB\times \R^d\times \R^q$ where $v_i(\theta):=\E[|\theta_n-\theta^*||\theta_i=\theta]$.
Under the assumptions of Theorem \ref{THM_RM}, the following control holds. %concerning the Lipschitz constant of $f_i^\gamma $.
\begin{PROP}[Controls of the Lipschitz constants]
\label{CTR_FIG}
%Assume \A{HL} and \A{HUA} are in force. Then, 
For all $i\in \leftB 1,n\rightB$, the function $f_i^\gamma $ satisfies:
$$
[f_i^\gamma]_1:=\sup_{\theta\in \R^d, y\neq y'}\frac{|f_i^\gamma(\theta,y)-f_i^\gamma(\theta,y')|}{|y-y'|}\le \left(\Pi_{n} \Pi_{i}^{-1} \right)^{1/2} \gamma_i [H]_1 .
$$

\noindent where $\Pi_{n} := \prod_{k=0}^{n-1} \left( 1 - 2 \underline{\lambda} \gamma_{k+1} + [H]^{2}_{1} \gamma^{2}_{k+1}\right)$, $n\geq1$.
\end{PROP}
The proof is postponed to  Section \ref{CTR_LIP_AS}.

\noindent %Hence, we are now in position to apply Proposition \ref{HERBST_COND} to derive the non-asymptotic bound for stochastic approximation algorithms. 
Set $\scD_{i}=f_i^{\gamma}(\theta_{i-1},Y_{i})- \E[ \left. f_i^{\gamma}(\theta_{i-1},Y_{i}) \right| \mathcal{F}_{i-1}]$. Recalling that the random variables $(Y_i)_{i\in \N^*}$ satisfy \eqref{CONC_GAUSS},  we obtain from Proposition \ref{HERBST_COND} that for all $r \geq 0$:
\begin{eqnarray*}
\P\left(\left|\theta_{N} - \theta^{*}\right| \geq r + \delta_{N}\right) =\P\left( \left|\theta_{N} - \theta^{*}\right| -\E[\left|\theta_{N} - \theta^{*}\right|]
\geq r \right)\\
\leq \exp\left(- \frac{r^{2}}{\alpha [H]^{2}_{1}\Pi_{N}\sum_{k=1}^{N}\gamma^{2}_{k}/\Pi_{k}}\right).
\end{eqnarray*}

\medskip

Contrary to the result concerning the Euler scheme, a bias appears in the non-asymptotic bound for the stochastic approximation algorithm. Consequently, it is crucial to have a control on it. At step $n$ of the algorithm, it is equal to $\delta_{n}:=\E[\left|\theta_{n}-\theta^{*}\right|]$. Under the current assumptions \A{HL} of \textit{Lipschitz continuity} of $H$ and \A{HUA} of \textit{uniform attractivity}, we have the following proposition.
\begin{PROP}[Control of the bias]
\label{CTR_BIAS}
For all $n \ge 1$, we have 
\begin{eqnarray*}
\delta_{n} \leq \exp\left(-\underline{\lambda} \Gamma_{n}\right) \left|\theta_0 - \theta^{*}\right|
+ [H]_1\sigma_{Y} \left(\sum_{k=0}^{n-1} e^{-2\underline{\lambda} (\Gamma_{n} - \Gamma_{k+1})} \gamma^{2}_{k+1} \right)^{\frac12},
\end{eqnarray*}

\noindent where $\Gamma_{n} :=\sum_{k=1}^{n} \gamma_{k}$, $\sigma_{Y} := \E\left[F^{2}(Y)\right]^{1/2}<+\infty$, with $F:y \mapsto \E\left[|y-Y|\right]$.
\end{PROP}

\noindent \textit{Proof.}
With the notations of Section \ref{RM_SEC}, we define for all $n\ge 1, \ \Delta M_n:= h(\theta_{n-1})- H(\theta_{n-1},Y_n)= \E[\left. H(\theta_{n-1},Y_n) \right| \mathcal{F}_{n-1}]- H(\theta_{n-1},Y_n)$. Recalling that $(Y_i)_{i\in \N^*} $ is a sequence of i.i.d. random variables we have that $(\Delta M_n)_{n\ge 1} $
is a sequence of martingale increments w.r.t. the natural filtration $\left(\F_n:=\sigma(Y_i, i \leq n), n\ge 1 \right) $.

From the dynamics \eqref{RM}, write now for all $n\in \N$,
\begin{eqnarray*}
z_{n+1} & := & \theta_{n+1}-\theta^*=\theta_n-\theta^*- \gamma_{n+1}\left\{h(\theta_n)-\Delta M_{n+1}\right\}\\
 & = & \theta_{n}-\theta^*-\gamma_{n+1}  \int_0^1d\lambda Dh(\theta^*+\lambda  (\theta_n-\theta^*) )(\theta_n-\theta^*) +\gamma_{n+1}\Delta M_{n+1},
\end{eqnarray*}

\noindent where we used that $h(\theta^*)=0 $ for the last equality. Setting $J_n:=\int_0^1d\lambda Dh(\theta^*+\lambda (\theta_n-\theta^*) )$, we obtain
\begin{eqnarray*}
z_{n+1} & = & (I-\gamma_{n+1} J_n)z_n + \gamma_{n+1} \Delta M_{n+1}.
\end{eqnarray*}

\noindent Take now the square of the $L^2$-norm in the previous equality. Recalling that  $\Delta M_{n+1}$ is a martingale increment, we derive:
\begin{align*}
\E[|z_{n+1} |^2] & = \E[|(I-\gamma_{n+1} J_n)z_n|^2] + \gamma^2_{n+1} \E[|\Delta M_{n+1}|^2] \\
&  \le \exp(-2 \underline{\lambda} \gamma_{n+1})\E[ |z_n|^2] +  \gamma^2_{n+1} [H]^2_1 \sigma^2_Y.
\end{align*}

\noindent For the last inequality we used (exploiting assumption \A{HUA}, uniform attractivity of the Jacobian matrix of $h$) $||I-\gamma_{n+1} J_n|| \le \exp(-\underline{\lambda} \gamma_{n+1})$, $\|.\| $ standing for the matrix norm on $\R^d\otimes \R^d $, and the inequality $\E[|\Delta M_{n+1}|^2]\le [H]_1^2\sigma_Y^2 $ which follows from \A{HL}.

A direct induction yields for all $n\ge 1$:
\begin{eqnarray*}
\E[|z_n|^2]\le \exp\left(-2 \underline{\lambda }\Gamma_n \right)|z_0|^2+[H]_1^2\sigma_Y^2\left(\sum_{k=0}^{n-1} e^{-2\underline{\lambda}(\Gamma_n-\Gamma_{k+1})} \gamma_{k+1}^2\right),
\end{eqnarray*}
which completes the proof.
\appendix

\section{Technical results}
\subsection{Proof of Proposition \ref{PROP_LIP_EUL}}
\label{CTR_LIP}
The proof follows from usual stochastic analysis arguments that we now recall for the sake of completeness. 
For $i=N$ we directly get from the definition of $f_N^\Delta(x,y)$ that $[f_N^\Delta]_1\le [f]_1|\sigma|_\infty $.

For $i\in \leftB 1,N-1\rightB $, define
for $y\ne y' $ the quantity 
\begin{eqnarray*}
D_{t_i}^{\Delta}(T,x,y,y')&:=&\sup_{j\in \leftB i,N\rightB}\frac{|X_{t_j}^{\Delta,t_i,G_{i-1}^\Delta( x,y)} -X_{t_j}^{\Delta,t_i,G_{i-1}^\Delta(x,y')}|}{|y-y'|},\\
  \forall z\in \R^d,\ G_{i-1}^\Delta(x,z)&:=&b(t_{i-1},x)\Delta+\sigma(t_{i-1},x) z.
\end{eqnarray*}
%%%% Il y a une petite etape de localisation avec Fatou a faire pour l'aspect fini, ce qui precedait n'etait pas coherent.

 Write now:
 \begin{eqnarray*}
 &&D_{t_i}^\Delta(T,x,y,y')\\
 &\le& \left \{ \frac{|G_{i-1}^\Delta(x,y)-G_{i-1}^\Delta(x,y')|}{|y-y'|}+\Delta [b]_1\left(\bsum{k=i}^{N-1} \frac{|X_{t_k}^{\Delta, t_i,G_{i-1}^\Delta(x,y)}-X_{t_k}^{\Delta, t_i,G_{i-1}^\Delta(x,y')}|}{|y-y'|}\right) \right.\\
 &&+\left. \sup_{j\in \leftB i,N\rightB}\left|\bsum{k=i}^{j} \left (\frac{\sigma(t_k,X_{t_k}^{\Delta,t_i,G_{i-1}^{\Delta}(x,y) })-\sigma(t_k,X_{t_k}^{\Delta,t_i,G_{i-1}^{\Delta}(x,y')} )}{|y-y'|} \right) \sqrt \Delta Y_{k+1} \right|\right \}.
 \end{eqnarray*}
One can easily show that for all $(z,j)\in \R^q\times\leftB 0,N\rightB ,\ X_{t_j}^{\Delta,t_i,G_{i-1}^\Delta( x,z)}\in L^1(\P)$. Hence, taking the expectation and using the the Burkholder-Davis-Gundy inequality (see e.g. Chapter 7, \S 3 in Shiryaev \cite{shiryaev:96}) we obtain:
 \begin{eqnarray}
 &&\E[D_{t_i}^\Delta(T,x,y,y')]\nonumber\\
 &\le&\left\{ |\sigma|_\infty+ [b]_1\Delta \sum_{k=i}^{N-1}  \E\left[\frac{|X_{t_k}^{\Delta, t_i,G_{i-1}^\Delta(x,y)}-X_{t_k}^{\Delta, t_i,G_{i-1}^\Delta(x,y')}|}{|y-y'|}\right]\right. \nonumber \\
 &&
+   \left. c \E\left[\left( \sum_{k=i}^{N-1} \left| \frac{\sigma(t_k,X_{t_k}^{\Delta,t_i,G_{i-1}^{\Delta}(x,y) })-\sigma(t_k,X_{t_k}^{\Delta,t_i,G_{i-1}^{\Delta}(x,y')} )}{|y-y'|}\sqrt \Delta Y_{k+1}\right|^2\right)^{1/2}\right]   \right\}.\nonumber \\
&\le & |\sigma|_\infty+[b]_1\Delta \bsum{k=i}^{N-1}\E[(D_{t_i}^\Delta(t_k,x,y,y')]\nonumber\\
&&+c[\sigma]_1\E\left[ \left(\Delta \bsum{k=i}^{N-1} \frac{|X_{t_k}^{\Delta,t_i,G_{i-1}^\Delta(x,y)}-X_{t_k}^{\Delta,t_i,G_{i-1}^\Delta(x,y)}|^2  }{|y-y'|^2}|Y_{k+1}|^2\right)^{1/2}\right] ,\nonumber\\
&& c:=c(q).\label{PRELIM_GRONW}
\end{eqnarray}
Observe now that 
\begin{eqnarray*}
\E\left[ \left(\Delta \bsum{k=i}^{N-1} \frac{|X_{t_k}^{\Delta,t_i,G_{i-1}^\Delta(x,y)}-X_{t_k}^{\Delta,t_i,G_{i-1}^\Delta(x,y)}|^2  }{|y-y'|^2}|Y_{k+1}|^2\right)^{1/2} \right]\\
\le \E\left[D_{t_i}^\Delta(T,x,y,y')^{1/2}\left(\Delta \bsum{k=i}^{N-1} D_{t_i}^\Delta(t_k,x,y,y') |Y_{k+1}|^2\right)^{1/2}\right] \\
\le \eta \E[D_{t_i}^\Delta(T,x,y,y')]+\eta^{-1} \bsum{k=i}^{N-1} \E[D_{t_i}^\Delta(t_k,x,y,y')]\E[|Y_1|^2],\forall \eta\in (0,1),
\end{eqnarray*}
which plugged into \eqref{PRELIM_GRONW}  yields thanks to the Gronwall Lemma 
\begin{eqnarray*}
(1-c[\sigma]_1\eta)\E[D_{t_i}^\Delta(T,x,y,y')]\le |\sigma|_\infty \exp\left(  \{[b]_1+ c [\sigma]_1\E[|Y_1|^2]\eta^{-1} \} (T-t_i)\right),\\
 \eta \in \left(0,\frac{1}{c[\sigma]_1}\wedge 1\right).
\end{eqnarray*}
Taking  $ \eta:=\frac{(c[\sigma]_1)^{-1}\wedge 1}2 $ we obtain
$$\E[D_{t_i}^\Delta(T,x,y,y')]\le 2 |\sigma|_\infty \exp\left(  \{[b]_1+ 2c \E[|Y_1|^2] [\sigma]_1 (1\vee c[\sigma]_1)  \} (T-t_i)\right),$$
which recalling $[f_i^\Delta]_1\le [f]_1 \sup_{y\neq y'}\E[|D_{t_i}^\Delta(T,x,y,y' )|]$ %and %$\E[|Y_1|^2]=q $ (indeed we assumed in Section \ref{DESCR_EUL} that $\E[Y_1Y_1^*]=I_q $) 
completes the proof up to a modification of $c$.
%%%% Il faut reporter ceci i.e. le fait qu'apparaisse un \E[|Y_1|^2] dans la constante, le c est maintenant juste la constante de BDG. Auparavant cette depedendance etait cachee dans le c

\subsection{Proof of Proposition \ref{CTR_FIG}}
\label{CTR_LIP_AS}
%\noindent \textit{Proof} 
From the definitions in Section \ref{RM_SEC}, it suffices to control the difference $\E[|\theta_n^{\theta,i}-\theta_n^{\theta',i}|]$, that is the sensitivity of the algorithm w.r.t. the starting point at time $i\in \leftB 1,n\rightB$. Write for all $j\in \leftB i,n-1\rightB $: 
\begin{eqnarray*}
|\theta_{j+1}^{\theta,i}-\theta_{j+1}^{\theta',i}|^2 & = & |\theta_{j}^{\theta,i}-\theta_{j}^{\theta',i}-\gamma_{j+1}\left\{H(\theta_{j}^{\theta,i},Y_{j+1})-H(\theta_{j}^{\theta',i},Y_{j+1}) \right\}|^2\\
& = & |\theta_{j}^{\theta,i}-\theta_{j}^{\theta',i}|^2 - 2 \gamma_{j+1} \langle \theta_{j}^{\theta,i}-\theta_{j}^{\theta',i},H(\theta_{j}^{\theta,i},Y_{j+1})-H(\theta_{j}^{\theta',i},Y_{j+1})\rangle\\
&& + \gamma_{j+1}^2|H(\theta_{j}^{\theta,i},Y_{j+1})-H(\theta_{j}^{\theta',i},Y_{j+1})|^2 \\
& = & |\theta_{j}^{\theta,i}-\theta_{j}^{\theta',i}|^2 - 2 \gamma_{j+1} \langle \theta_{j}^{\theta,i}-\theta_{j}^{\theta',i},h(\theta_{j}^{\theta,i})-h(\theta_{j}^{\theta',i})\rangle \\
&& - 2 \gamma_{j+1} \langle \theta_{j}^{\theta,i}-\theta_{j}^{\theta',i}, \Delta M^{\theta,i}_{j+1} - \Delta M^{\theta',i}_{j+1}\rangle \\
&&+ \gamma_{j+1}^2|H(\theta_{j}^{\theta,i},Y_{j+1})-H(\theta_{j}^{\theta',i},Y_{j+1})|^2,
\end{eqnarray*}
%%% Rq j'ai change le + en - pour la somme des increments de martingale.

\noindent where we introduced the martingale increments $\Delta M^{\theta,i}_{j+1} = H(\theta_{j}^{\theta,i},Y_{j+1}) - h(\theta_{j}^{\theta,i})$ and $\Delta M^{\theta',i}_{j+1} = H(\theta_{j}^{\theta',i},Y_{j+1}) - h(\theta_{j}^{\theta',i})$, $j\geq0$ in the last equality. Now, using \A{HL} and \A{HUA} yields:
\begin{eqnarray*}
|\theta_{j+1}^{\theta,i}-\theta_{j+1}^{\theta',i}|^2 & \leq & |\theta_{j}^{\theta,i}-\theta_{j}^{\theta',i}|^2 \left(1 -2 \underline{\lambda} \gamma_{j+1} + [H]^{2}_{1} \gamma_{j+1}^{2}\right) \\
&&- 2 \gamma_{j+1} \langle \theta_{j}^{\theta,i}-\theta_{j}^{\theta',i}, \Delta M^{\theta,i}_{j+1} - \Delta M^{\theta',i}_{j+1}\rangle,   
\end{eqnarray*}

\noindent and, by induction on $j$, we easily obtain:
\begin{eqnarray}
|\theta_n^{\theta,i}-\theta_n^{\theta',i}|^2  \leq  |\theta-\theta'|^2 \prod_{j=i}^{n-1} \left(1 -2 \underline{\lambda} \gamma_{j+1} + [H]^{2}_{1} \gamma^{2}_{j+1}\right) \nonumber \\ 
 - 2 \left(\prod_{j=i}^{n-1} \left(1 -2 \underline{\lambda} \gamma_{j+1} + [H]^{2}_{1} \gamma^{2}_{j+1}\right) \right)\sum_{j=i}^{n-1} \tilde{\gamma}_{j+1} \langle \theta_{j}^{\theta,i}-\theta_{j}^{\theta',i}, \Delta M^{\theta,i}_{j+1} - \Delta M^{\theta',i}_{j+1}\rangle \label{PREAL_GRON_ALGO_STO}
\end{eqnarray}

\noindent where $\tilde{\gamma}_{j+1} := \gamma_{j+1} / \prod_{k=i}^{j} \left(1 -2 \underline{\lambda} \gamma_{k+1} + [H]^{2}_{1} \gamma^{2}_{k+1}\right)$. Taking the expectation in \eqref{PREAL_GRON_ALGO_STO}, we derive:
$$
\frac{\E[|\theta_{n}^{i,\theta}-\theta_{n}^{i,\theta'}|^2]}{|\theta-\theta'|^2} \le \prod_{j=i}^{n-1} \left(1 -2 \underline{\lambda} \gamma_{j+1} + [H]^{2}_{1} \gamma^{2}_{j+1}\right).
$$

\noindent Now:
\begin{eqnarray*}
|f_i^\gamma(\theta,y)-f_i^\gamma(\theta,y')| & = & |\E[|\theta_n^{i,\theta-\gamma_i H(\theta,y)}-\theta^{*}|] -\E[|\theta_n^{i,\theta-\gamma_i H(\theta,y')}-\theta^{*}|]|\\
&\le &\E[|\theta_n^{i,\theta-\gamma_i H(\theta,y)}-\theta_n^{i,\theta-\gamma_i H(\theta,y')}|]\\
& \le & \left(\prod_{j=i}^{n-1} \left(1 -2 \underline{\lambda} \gamma_{j+1} + [H]^{2}_{1} \gamma^{2}_{j+1}\right)\right)^{1/2} \gamma_i [H]_1 |y-y'| \\
& = & \left(\Pi_{n} \Pi_{i}^{-1}\right)^{1/2} \gamma_i [H]_1 |y-y'|,
\end{eqnarray*}

\noindent which completes the proof of the Proposition.

%%%%%%%%%%%%%%%%%%%%%%%%%%%%%%%%%%%%%%%%%%%%%%%%%%%%%%%%%%%%%%%%%%%
%%                                                               %%
%% Use the two commands below for producing your bibliography    %%
%% with bibtex, then comment again the commands and include the  %%
%% content of the .bbl file in this file below the commands.     %%
%%                                                               %%
%%%%%%%%%%%%%%%%%%%%%%%%%%%%%%%%%%%%%%%%%%%%%%%%%%%%%%%%%%%%%%%%%%%

%\bibliographystyle{amsplain}
%\bibliography{bibli}

\providecommand{\bysame}{\leavevmode\hbox to3em{\hrulefill}\thinspace}
\providecommand{\MR}{\relax\ifhmode\unskip\space\fi MR }
\providecommand{\MRhref}[2]{%
  \href{http://www.ams.org/mathscinet-getitem?mr=#1}{#2}
}
\providecommand{\href}[2]{#2}

% add below the content of your .bbl file produced by bibtex.

%\begin{thebibliography}{99}
%
%\bibitem{doob} Doob, J. L.: Heuristic approach to the Kolmogorov-Smirnov
%  theorems. \emph{Ann. Math. Statistics} \textbf{20}, (1949), 393--403.
%  \MR{0030732}
%
%\bibitem{gnekol} Gnedenko, B. V. and Kolmogorov, A. N.: Limit distributions for
%  sums of independent random variables. Translated and annotated by K. L.
%  Chung. With an Appendix by J. L. Doob. \emph{Addison-Wesley}, Cambridge,
%  1954. ix+264 pp. \MR{0062975}
%
%\bibitem{ito} It\^o, K.: Multiple Wiener integral. \emph{J. Math. Soc. Japan}
%  \textbf{3}, (1951), 157--169. \MR{0044064}
%
%\bibitem{levy} L\'evy, P.: Sur certains processus stochastiques homog\`enes.
%  \emph{Compositio Math.} \textbf{7}, (1939), 283--339. \MR{0000919}
%
%\bibitem{grisha} Perelman, G.: The entropy formula for the Ricci flow and its
%  geometric applications, \ARXIV{math.DG/0211159}
%
%\bibitem{smisch} Smirnov, S. and Schramm, O.: On the scaling limits of planar
%  percolation, \ARXIV{1101.5820}
%
%\end{thebibliography}

%%%%%%%%%%%%%%%%%%%%%%%%%%%%%%%%%%%%%%%%%%%%%%%%%%%%%%%%%%%%%%%%%%%
%%                                                               %%
%% You may add acknowledgments (optional).                       %%
%%                                                               %%
%%%%%%%%%%%%%%%%%%%%%%%%%%%%%%%%%%%%%%%%%%%%%%%%%%%%%%%%%%%%%%%%%%%

\section*{Acknowledgments}We would like to thank Florent Malrieu for fruitful discussions and two
  anonymous referees for their careful reading and comments. We thank Bernard Bercu who pointed out an error in the control of the bias in the previously published version of this work.

\end{document}